\documentclass[12pt]{article}
\usepackage{amsmath,amsfonts,amsthm,amscd,amssymb,latexsym,comment}
\usepackage{enumerate,graphicx,psfrag,subfigure,changebar,oldgerm}
\title{Parabolic and Quasiparabolic Subgroups of Free Partially Commutative
Groups}
\author{ \textsf{Andrew J. Duncan}
 \and \textsf{Ilya V. Kazachkov}
 \and \textsf{Vladimir N. Remeslennikov}}

\def\nul{\emptyset }
\def\D{\Delta }

\def\G{\Gamma }

\def\a{\alpha }

\def\fC{{\textswab C}}
\newcommand{\gd}{\operatorname{gd}}
\renewcommand{\gcd}{\operatorname{gd}}
\renewcommand{\lg}{{l}}
\newtheorem{thm}{Theorem}[section]
\newtheorem{lem}[thm]{Lemma}
\newtheorem{cor}[thm]{Corollary}
\newtheorem{prop}[thm]{Proposition}

\newtheorem{defn}[thm]{Definition}
\newtheorem{exam}[thm]{Example}

\numberwithin{equation}{section} \numberwithin{figure}{section}

\newcommand{\cB}{\mathcal{B}}

\newcommand{\cO}{\mathcal{O}}
\newcommand{\cP}{\mathcal{P}}

\newcommand{\az}{\mathop{{\ensuremath{\alpha}}}}
\newcommand{\la}{\langle}
\newcommand{\ra}{\rangle}

\newcommand{\rank}{\operatorname{rank}}
\newcommand{\cl}{\operatorname{cl}}

\newcommand{\CS}{\operatorname{CS}}

\newcommand{\ov}[1]{\overline{#1}}
\newcommand{\bs}{\setminus}
\newcommand{\be}{\begin{enumerate}}
\newcommand{\ee}{\end{enumerate}}

\newcommand{\ilya}[1]{\marginpar{{\LARGE $\mathbf{\Downarrow IVK}$}}\textbf{ #1}}
\begin{document}

\maketitle

\begin{abstract}
Let $\Gamma$ be a finite graph and $G$ be the corresponding
free partially commutative group. In this paper we study subgroups 
generated by vertices of the graph $\Gamma$, which we call  canonical parabolic subgroups.
A natural extension of the definition 
leads to canonical quasiparabolic subgroups. It is shown that
the centralisers of subsets of $G$ are the conjugates of
canonical quasiparabolic centralisers satisfying certain graph theoretic
conditions.
\end{abstract}

 \tableofcontents
\section{Preliminaries}
In this section we give a brief overview of some definitions and
results from \cite{EKR,DKR3}.
We begin with the basic notions of the theory of free partially
commutative groups.
Let $\G$ be a finite, undirected, simple graph. Let
$X=V(\G)=\{x_1,\dots, x_n\}$ be the set of vertices of $\G$ and let
$F(X)$ be the free group on $X$. Let
\[
R=\{[x_i,x_j]\in F(X)\mid x_i,x_j\in X \textrm{ and there is an edge of }
\G \textrm{ joining }
x_i \textrm{ to } x_j \}.
\]
We define the {\em partially commutative group} with ({\em commutation}) 
{\em graph}
$\G$ to be the group $G(\G)$ with presentation $ \left< X\mid
R\right>$. When the underlying graph is clear from the context we
write simply $G$.

Denote by $\lg(g)$ the minimum of the lengths words that represents the
element $g$. If $w$ is a word representing $g$ and $w$ has length
$\lg(g)$ we call $w$ a {\em minimal form} for $g$. 
When the meaning is clear we shall say that $w$ is a 
minimal element of $G$ when we mean that $w$ is a minimal form of an
element of $G$.
We say that   $w \in G$
 is {\em cyclically minimal} if and only if
$$
\lg(g^{-1}wg) \ge \lg(w)
$$
for every $g \in G$. We write $u\circ w$ to express the fact that
$\lg(uw)=\lg(u)+\lg(w)$, where $u,w\in G$.
We will need the notions of a divisor and the greatest
divisor of a word $w$ with respect to a subset $Y\subseteq X$, defined
in \cite{EKR}.
Let $u$ and $w$ be elements of $G$. We say that $u$ is a {\em left }
({\em right}) {\em divisor} of $w$ 
if there exists $v\in G$ such that $w= u
\circ v$ ($w=  v\circ u$).
We order the set of all left (right) divisors of a word $w$ as
follows. We say that $u_2$ is greater than $u_1$ if and only if
$u_1$ left (right) divides $u_2$.
It is shown in \cite{EKR} that, for any $w\in G$ and $Y\subseteq
X$, there exists a  unique maximal left  divisor  of $w$ 
which belongs to the subgroup $G(Y) <G$ which is called 
the {\em greatest left divisor} $\gd^{l}_Y(w)$ of $w$ in $Y$.
  The {\em greatest right divisor} of $w$ in $Y$ is defined analogously.
 We omit the indices when no ambiguity occurs.

The {\em non-commutation} graph of 
the partially commutative  group $G(\G)$ is the graph $\Delta$,
dual to $\Gamma$, with vertex set $V(\Delta)=X$ and 
an edge connecting $x_i$ and $x_j$ if and only if $\left[x_i, x_j
\right] \ne 1$. The graph $\D$ is a union of its connected
components $\D_1, \ldots , \D_k$ and  words that depend on
letters from distinct components commute. For any graph $\G$, 
if $S$ is a subset
of $V(\G)$ we shall write $\G(S)$ for the full subgraph of $\G$ with
vertices $S$. Now, if
the vertex set of 
$\D_k$ is $I_k$ and $\G_k=\G(I_k)$ 
then 
$G=
G(\G_1) \times \cdots \times G(\G_k)$.
For $g \in G$ let $\az(g)$ be the set of elements $x$ of $X$ such
that $x^{\pm 1}$ occurs in a minimal word $w$ representing $g$.
It is shown in \cite{EKR} that $\a(g)$ is well-defined.
Now suppose that  the full subgraph $\Delta (\az(w))$ 
of $\Delta$ with vertices $\a(w)$ has connected components
$\D_1,\ldots, \D_l$ and let the vertex set of $\D_j$ be $I_j$.
Then, since $[I_j,I_k]=1$, 
we can split $w$ into the product of commuting words,
$w=w_1\circ \cdots \circ w_l$, where $w_j\in G(\G(I_j))$, so $[w_j,w_k]=1$
for all $j,k$.
If $w$ is cyclically minimal then
we call this expression for $w$ a
 {\em block decomposition} of $w$ and
say $w_j$ a {\em block} of $w$, for $j=1,\ldots ,l$. Thus $w$ itself is
 a block if and only if  $\D(\a(w))$ is connected. In general let 
$v$ be an element of $G$ which is not necessarily cyclically minimal. We may write 
$v=u^{-1}\circ w \circ u$, where $w$ is cyclically minimal and then
$w$ has a block decomposition $w=w_1\cdots w_l$, say. Then
$w_j^u=u^{-1}\circ w_j\circ u$ and we call the expression $v=w_1^u \cdots w_l^u$ 
the  {\em block decomposition}
of $v$ and say that $w_j^u$ is a {\em block}
of $v$, for $j=1,\ldots , l$. Note that this definition is slightly different
from that given in \cite{EKR}.

Let $Y$ and $Z$ be subsets of $X$. As in \cite{DKR3} 
we define the {\em orthogonal complement}
of $Y$ {\em in} $Z$ to be
\[\cO^Z(Y)=\{u\in Z|d(u,y)\le 1, \textrm{ for all } y\in Y\}.\]
By convention we set $\cO^Z(\nul)=Z$.
If $Z=X$ we call $\cO^X(Y)$ the orthogonal complement of $Y$, and if
no ambiguity arises then we write $Y^\perp$ instead
of $\cO^X(Y)$ and $x^{\perp}$ for $\{x\}^\perp$. 
Let 
$\CS(\Gamma)$ be the set of all subsets $Z$ of $X$ of the form
$Y^\perp$ for some $Y\subseteq X$. The set $\CS(\G)$ is shown in \cite{DKR3} 
to
 be a lattice, the {\em lattice of
closed sets} of $\Gamma$. 

The {\em centraliser} of a subset $S$ of  $G$  is 
\[C(S)=C_G(S)=\{g\in G: gs=sg,
\textrm{ for all } s\in S\}.\] 
The set $\fC(G)$ of centralisers of a group is a lattice. 
An element $g \in G$ is called a {\em root element} 
if $g$ is not a proper
power of any element of $G$. If $h=g^n$, where $g$ is a root element
and $n\ge 1$, then $g$ is said to be a {\em root} of $h$.
As shown in \cite{DK} every element of the partially commutative group
 $G$ has a unique root, which
we denote $r(g)$. 
If $w\in G$ define $A(w)=\langle Y\rangle=G(Y)$, where
$Y=\az(w)^\perp \bs \az(w)$. 
Let $w$ be a cyclically
minimal element of $G$ with block decomposition $w=w_1\cdots w_k$
and let $v_i=r(w_i)$. Then, from \cite[Theorem 3.10]{DK},
\begin{equation}\label{eq:centraliser}
C(w)=\langle v_1\rangle \times \cdots \times \langle v_k
\rangle\times A(w).
\end{equation}

We shall use \cite[Corollary 2.4]{DKR2} several times in what follows, 
so for ease of reference we state it here: first
recalling the necessary notation.
It follows from
\cite[Lemma 2.2]{DKR2} that if $g$ is a cyclically minimal element of $G$ 
and $g=u\circ v$ then $vu$ is cyclically minimal.
%
For a cyclically minimal element $g\in G$ we define 
$\tilde{g}=\{h\in G|h=vu, \textrm{
for some } u,v \textrm{ such that } g=u\circ v\}$. (We allow $u=1$, $v=g$
so that $g\in \tilde{g}$.)

\begin{lem}{\bf \cite[Corollary 2.4]{DKR2}}\label{cor:conj}
Let $w,g$ be (minimal forms of) elements of $G$ and
$w=u^{-1}\circ v\circ u$, where $v$ is cyclically minimal. Then
there exist minimal forms $a$, $b$, $c$, $d_1$, $d_2$ and
$e$ such that $g=a\circ b\circ c\circ d_2$, $u=d_1\circ a^{-1}$, $d=d_1\circ d_2$,
$w^g=d^{-1}\circ e\circ d$, $\tilde{e}=\tilde{v}$, $e=v^b$, $\a(b)\subseteq \a(v)$ and
$[\a(b\circ c),\a(d_1)]=[\a(c),\a(v)]=1$.
\end{lem}

\begin{figure}\label{fig:vk}
\begin{center}
\psfrag{a}{$a$}
\psfrag{b}{$b$}
\psfrag{c}{$c$}
\psfrag{d1}{$d_1$}
\psfrag{d2}{$d_2$}
\psfrag{f}{$f$}
\psfrag{0}{$0$}
\psfrag{A}{$A$}
\psfrag{B}{$B$}
\includegraphics[scale=0.4]{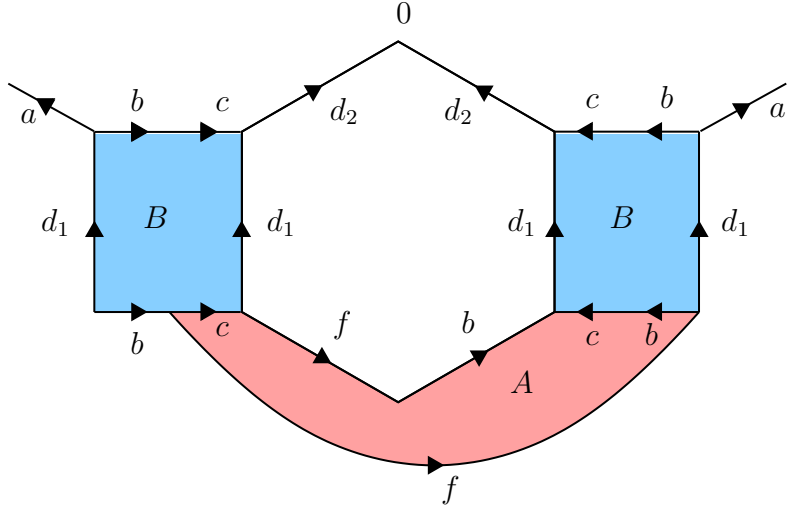}
\end{center}
\caption{A Van Kampen diagram for Lemma \ref{cor:conj}}
\end{figure}

Figure \ref{fig:vk} expresses the conclusion
of  Lemma \ref{cor:conj} as a Van Kampen diagram. In this diagram
we have assumed that $v=b\circ f$ and so $e=f\circ b$. 
The regions labelled $B$ are tessellated using relators corresponding
to the relation $[\a(b\circ c),\a(d_1)]=1$ and the region labelled $A$ 
with relators corresponding to $[\a(c),\a(v)]=1$. 
Reading
anticlockwise from the vertex labelled $0$  the boundary
label of the exterior region is $g^{-1}wg$ and the label of the interior
region (not labelled $A$ or $B$) is $e^d$.

\section{Parabolic Subgroups}

\subsection{Parabolic and Block-Homogeneous Subgroups}

As usual let $\G$ be a graph with vertices $X$ and $G=G(\G)$. If $Y$
is a subset of $X$ denote by $\G(Y)$ the full subgraph of $\G$ with
vertices $Y$.
 Then $G(\G(Y))$ is the
free partially commutative group with graph $\G(Y)$. In \cite{EKR}
it is shown that $G(\G(Y))$ is the subgroup $\la Y \ra$ of $G(\G)$
generated by $Y$. We call $G(\G(Y))$ a {\em canonical parabolic}
subgroup of $G(\G)$ and, when no ambiguity arises, denote it $G(Y)$.
The elements of $Y$ are termed the {\em canonical} generators of
$G(Y)$.

\begin{defn}
A subgroup $P$ of  $G$ is called {\em parabolic} if it is conjugate
to a canonical parabolic subgroup $G(Y)$ for some $Y\subseteq X$.
The {\em rank} of $P$ is the cardinality $|Y|$ and $Y$ is called a
set of {\em canonical generators} for $P$.
\end{defn}

To see that the definition of rank of a parabolic subgroup is well
defined  note that if $Y,Z\subseteq X$ and $G(Y)=G(Z)^g$, for some
$g\in G$, then we have $y=g^{-1}w_yg$, for some $w_y\in G(Z)$, for all $y\in Y$.
It follows, from \cite[Lemma 2.5]{EKR}, by counting the exponent 
sums of letters
in a geodesic word representing $g^{-1}wg$, 
that $y\in \a(w_y)$, so $y\in Z$. Hence
$Y\subseteq Z$ and similarly $Z\subseteq Y$ so $Y=Z$. 

\begin{defn}
  A subgroup $H$ is called {\em block-homogeneous} if, for all $h\in H$,
if $h$ has block decomposition $h =w_1w_2 \ldots w_k$ then
$w_i \in H$, for $i=1,\ldots, k$.
\end{defn}

\begin{lem} \label{lem:4.1}
  An intersection of block-homogeneous subgroups is again a
block-homogeneous subgroup. 
If $H$ is block-homogeneous and $g\in G$ then $H^g$ is block-homogeneous.
In particular
parabolic subgroups are block-homogeneous.
\end{lem}
\begin{proof}
The first statement follows 
  directly from the definition. 
Let $H$ be block-homogeneous and $g\in G$ and let $w^g\in H^g$, where 
$w\in H$. Write $w=u^{-1}\circ v\circ u$, where $v$ is cyclically reduced and
has block-decomposition $v=v_1\cdots v_k$. Then the blocks of $w$ are 
$v_j^u$, so $v_j^u\in H$, for $j=1,\ldots ,k$. From Lemma \ref{cor:conj} there
exist $a, b, c, d_1, d_2, e$ such that 
$g=a\circ b\circ c\circ d_2$, $u=d_1\circ a^{-1}$, $d=d_1\circ d_2$,
$w^g=d^{-1}\circ e\circ d$, $\tilde{e}=\tilde{v}$, $e=v^b$, $\a(b)\subseteq \a(v)$ and
$[\a(b\circ c),\a(d_1)]=[\a(c),\a(v)]=1$. 
As $\tilde{e}=\tilde{v}$ it follows that
$\D(\a(e))=\D(\a(v))$ so $e$ has block-decomposition $e=e_1\cdots e_k$, 
where $\tilde{e_j}=\tilde{v_j}$. 
Therefore $w^g$ has block-decomposition $w^g=e^d=e_1^d\cdots e_k^d$.
Moreover
$e=v^b$ so $e_j=v_j^b$. Thus
\[e_j^d=e_j^{cd}=v_j^{bcd}=v_j^{d_1bcd_2}=v_j^{ug}\in H^g,
\]
which implies that $H^g$ is block-homogeneous. 
It follows from \cite[Lemma 2.5]{EKR} that any canonical parabolic subgroup
is block-homogeneous and this gives the  final statement.
\end{proof}

\subsection{Intersections of parabolic subgroups}
In this section we show that an intersection of parabolic subgroups is
again a parabolic subgroup. To begin with we establish 
some preliminary results.

\begin{lem}\label{lem:euclid}
Let $Y,Z\subseteq X$, let $w\in G(Y)$ and  let $g\in G(X)$ be such
that
\begin{equation*}
\gcd{}^l_{Y}(g)=\gcd{}^r_{Z}(g)=1.
\end{equation*}
\be
\item\label{it:euclid1} If $w^g\in G(Z)$ then
$g\in A(w)$ and $w\in G(Y)\cap G(Z)=G(Y\cap Z)$.
\item\label{it:euclid2}
 If $Y=Z$ and $g\in C(w)$ then $g\in A(w)$.
\ee
\end{lem}
\begin{proof}
For \ref{it:euclid1},
in the notation of Lemma \ref{cor:conj} we have $w=u^{-1}\circ v\circ u$,
$w^g=d_2^{-1}\circ d_1^{-1}\circ e\circ d_1\circ d_2$ and
$g=a\circ b\circ c\circ d_2$. Applying the conditions of this Lemma we obtain
$a=b=d_2=1$, $u=d_1$ and $e=v$ so $w^g=w$ and $g=c$. Moreover, from Lemma
\ref{cor:conj} again we obtain $[\a(g),\a(w)]=1$. If $x\in \a(w)\cap \a(g)$ this
means that $g=x\circ g^\prime$, with $x\in Y$, contradicting  the hypothesis on
$g$. Hence $ \a(w)\cap \a(g)=\emptyset$ and $g\in A(w)$.
Statement  \ref{it:euclid2}
follows from \ref{it:euclid1}.
\end{proof}
\begin{cor}\label{cor:parin}
Let $Y,Z\subseteq X$ and $g\in G$. If $G(Y)^g\subseteq G(Z)$ and
$\gcd^l_{Y^\perp}(g)=1$ then $Y\subseteq Z$ and $\a(g)\subseteq Z$.
\end{cor}
\begin{proof}
Assume first that $\gcd^l_Y(g)=\gcd_Z^r(g)=1$. Let $y\in Y$ and 
$w=y$ in Lemma \ref{lem:euclid};
so $y^g\in G(Z)$ implies that $g\in A(y)$ and $y\in Z$. This holds for
all $y\in Y$ so we have $Y\subseteq Z$ and $g\in A(Y)$. Hence, in this case, $g=1$.
Now suppose that $g=g_1\circ d$, where $g_1=\gcd_Y^l(g)$. Then $\gcd_{Y^\perp}^l(g)=1$
implies that $\gcd_{Y^\perp}^l(d)=1$. Now write $d=e\circ g_2$, where
$g_2=\gcd_Z^r(d)$. Then $G(Y)^g=G(Y)^d$ and $G(Y)^d=G(Z)$ implies $G(Y)^e=G(Z)$.
As $\gcd^l_{Y^\perp}(d)=1$ the same is true of $e$ and from the above we conclude that
$e=1$ and that $Y\subseteq Z$. Now $g=g_1\circ g_2$, where $\a(g_1)\subseteq Y
\subseteq Z$ and $\a(g_2)\subseteq Z$. Thus $\a(g)\subseteq Z$, as required.
\end{proof}

\begin{prop} \label{prop:intpar}
Let $P_1$ and $P_2$ be parabolic subgroups. Then $P=P_1\cap P_2$ is
a parabolic subgroup. If $P_1\nsubseteq P_2$ then the rank of $P$ is
strictly smaller than the rank of $P_1$.
\end{prop}
This lemma follows easily from the next more technical result.
\begin{lem}\label{lem:formpar}
Let $Y,Z\subset X$ and $g\in G$. Then
$$
G(Y)\cap G(Z)^g=G(Y\cap Z\cap T)^{g_2},
$$
where $g=g_1\circ d\circ g_2$, $\gcd_Z^l(d)=\gcd_Y^r(d)=1$, $g_1\in G(Z)$,
$g_2\in G(Y)$
and
$T=\alpha(d)^\perp$.
\end{lem}
\begin{proof}[Derivation of Proposition \ref{prop:intpar} 
from Lemma \ref{lem:formpar}]
Let $P_1=G(Y)^a$ and $P_2=G(Z)^b$, for some $a,b\in G$. Then
$P=\left(G(Y)\cap G(Z)^{ba^{-1}}\right)^a$, which is parabolic since
Lemma \ref{lem:formpar} implies that
$G(Y)\cap G(Z)^{ba^{-1}}$ is parabolic. 
Assume that the rank of $P$ is 
greater than or equal to the rank of $P_1$.
Let $g=ba^{-1}$. The rank of $P$ is equal to the rank of 
$G(Y)\cap G(Z)^g$ and, in
the notation of Lemma \ref{lem:formpar}, 
$G(Y)\cap G(Z)^g=G(Y\cap Z\cap T)^{g_2}$, where 
$g=g_1\circ d\circ g_2$, with $T=\a(d)^\perp$,  $g_2\in G(Y)$ and $g_1\in G(Z)$. Therefore
$Y\subseteq Y\cap Z\cap T$ which implies $Y\subseteq Z\cap T$.
Thus $G(Y)\subseteq G(Z\cap T)=G(Z\cap T)^d$ so 
$G(Y)=G(Y)^{g_2}\subseteq G(Z\cap T)^{dg_2}\subseteq G(Z)^{dg_2}=G(Z)^g$.
Hence $P_1\subseteq P_2$. 
\end{proof}
\begin{proof}[Proof of Lemma \ref{lem:formpar}] 
Let $g_1=\gd^l_Z(g)$ and write $g=g_1\circ g^\prime$. Let 
$g_2=\gd^r_Y(g^\prime)$ and write $g^\prime =d\circ g_2$. Then
$g_1, g_2$ and $d$ satisfy the conditions of the lemma. Set
$T=\a(d)^\perp$.
As $G(Y)\cap G(Z)^g=G(Y)\cap G(Z)^{dg_2}=G(Y)^{g_2}\cap G(Z)^{dg_2}=
\left(G(Y)\cap G(Z)^d\right)^{g_2}$ it suffices to show that
$G(Y)\cap G(Z)^d=G(Y\cap Z\cap T)$. 
If $d=1$ then $T=X$ and $G(Y)\cap G(Z)=G(Y\cap Z)$, so the result holds.
Assume then that $d\neq 1$. 
Let $p=w^d\in G(Y)\cap G(Z)^d$, with $w\in G(Z)$. Applying
Lemma \ref{lem:euclid} to $w^d\in G(Y)$ we have
$d\in A(w)$ and $w\in G(Z)\cap
G(Y)=G(Y\cap Z)$. Thus $w\in \a(d)^\perp=T$ and so
$w\in G(Y\cap Z\cap T)$.
This shows that $G(Y)\cap G(Z)^d\subseteq G(Y\cap Z\cap T)$ and as
the reverse inclusion follows easily the proof is complete.
\end{proof}

\begin{prop}
    The intersection of parabolic subgroups is a parabolic subgroup and
can be obtained as an intersection of a finite number of subgroups
from the initial set.
\end{prop}
\begin{proof}
  In the case of two parabolic subgroups the result follows from
Proposition \ref{prop:intpar}. Consequently, the statement also holds for a
finite family of parabolic subgroups. For the general case we
use Proposition \ref{prop:intpar} again, noting that a proper intersection
of two parabolic subgroups is a parabolic subgroup of lower rank, and
the result follows.
\end{proof}
As a consequence of this Proposition we obtain: given two parabolic
subgroups $P$ and $Q$ the intersection $R$ of all parabolic
subgroups containing $P$ and $Q$ is the unique minimal parabolic
subgroup containing both $P$ and $Q$. Define $P\vee Q=R$ and
$P\wedge Q=P\cap Q$.
\begin{cor}
The parabolic subgroups of $G$ with the operations $\vee$ and $\wedge$
above form a lattice.
\end{cor}

\subsection{The Lattice of Parabolic Centralisers}

Let $Z \subseteq X$. Then the subgroup $C_G(Z)^g$ is called a
parabolic centraliser.  As shown in \cite[Lemma 2.3]{DKR3} every
parabolic centraliser is a parabolic subgroup: in fact
$C_G(Z)^g=G(Z^\perp)^g$. The converse also holds as the following proposition shows.
\begin{prop}
A parabolic subgroup $G(Y)^g$, $Y\subseteq X$ is a centraliser if
and only if there exists $Z\subseteq X$ so that $Z^\perp=Y$. In this
case $G(Y)^g=C_G(Z^g)$.
\end{prop}
\begin{proof}
  It suffices to prove the proposition for $g=1$ only. Suppose that
there exists such a $Z$. It is then clear that $G(Y)\subseteq
C_G(Z)$. If $w\in G$, $w$ is a reduced word and
$\alpha(w)\nsubseteq Y$ then there exists $x\in \alpha (w)$ and
$z\in Z$ so that $[x,z]\ne 1$ and consequently, by \cite[Lemma 2.4]{EKR}, 
$[w,z]\ne 1$. Assume further that $G(Y)$ is a
centraliser of a set of elements $w_1, \ldots, w_k$ written in a
reduced form. Since for any $y\in Y$ holds $[y,w_i] = 1$ then again,
by \cite[Lemma 2.4]{EKR}, $[y,\alpha(w_i)]=1$. Denote
$Z=\bigcup\limits_{i=1}^k \alpha(w_i)$. We have $[y,z]=1$ for all
$z\in Z$ and consequently $Y\subseteq Z^\perp$. 
Conversely if 
$x\in Z^\perp$ 
then 
$x\in C_G(w_1, \ldots, w_k)$ 
so $x\in Y$.
\end{proof}

We now introduce the structure of a lattice on the set of all
parabolic centralisers.
As we have shown above the intersection of two parabolic subgroups
is a parabolic subgroup. So, we set $P_1 \wedge P_2=P_1\cap P_2$.
The most obvious way to define $P_1\vee P_2$ would be to set $P_1\vee P_2=
\la P_1, P_2\ra$. However, in this case $P_1\vee P_2$ is not
necessarily a centraliser, though it is a parabolic subgroup.
For any $S\subseteq G$ we define the $\overline S=\cap\{P: P$ is a parabolic centraliser
and $S\subseteq P\}$. Then $\ov S$ is the minimal parabolic centraliser containing
$S$; since intersections of centralisers are centralisers and intersections of parabolic
subgroups are parabolic subgroups. We now define
$P_1\vee P_2=
\overline{\la P_1, P_2\ra}$.

\section{Quasiparabolic subgroups}
\subsection{Preliminaries}

As before let $\Gamma$ be a finite graph with vertex set $X$ 
and  $G=G(\Gamma)$ be
the corresponding 
partially commutative group.

\begin{defn}
Let 
$w$ be a
cyclically minimal root element of $G$ with block decomposition
$w=w_1 \cdots w_k$ and let $Z$ be a subset of $X$ such that
$Z\subseteq \a(w)^\perp$. Then the 
subgroup $Q=Q(w,Z)=\langle w_1 \rangle \times \dots \times \langle
w_k\rangle \times G(Z)$ is called a 
{\em canonical
quasiparabolic subgroup} of $G$.
\end{defn}

Note that we may choose $w=1$ so that canonical parabolic
subgroups are canonical quasiparabolic subgroups.
Given a canonical quasiparabolic subgroup $Q(w,Z)$, with $w$ and $Z$
as above, we may reorder the $w_i$ so that $\lg(w_i)\ge 2$, for 
$i=1,\ldots ,s$ and $\lg(w_i)=1$, for $i=s+1,\ldots ,k$. Then
setting $w^\prime=w_1\cdots w_s$ and  
$Z^\prime=Z\cup \{w_{s+1},\ldots, w_k\}$ we have
$Z^\prime \subseteq \a(w)^\perp$ and $Q(w,Z)=Q(w^\prime,Z^\prime)$.
This prompts the following definition.

\begin{defn}
We say that a canonical quasiparabolic subgroup $Q=\langle w_1
\rangle\times \dots \times \langle w_k\rangle \times G(Z)$ is written
in standard form if $|\alpha(w_i)| \ge 2$, $i=1,\dots, k$, or
$w=1$. 
\end{defn}

There are two obvious advantages to the standard form which we record
in the following lemma.

\begin{lem}\label{lem:stdform}
The standard form of  a canonical quasiparabolic subgroup $Q$
is unique, up to reordering of blocks of $w$. If
$Q(w,Z)$ is the standard form of $Q$ then $Z\subseteq \a(w)^\perp\bs \a(w)$.
\end{lem}
\begin{proof}
That the standard form is unique follows from uniqueness of 
roots of elements in partially commutative groups. The second statement
follows directly from the definitions.
\end{proof}
\begin{defn}
A subgroup $H$ of \/ $G$ is called quasiparabolic if it is
conjugate to a canonical quasiparabolic subgroup.
\end{defn}

Let $H=Q^g$ be a quasiparabolic subgroup of $G$, where $Q$ is the
canonical quasiparabolic subgroup of $G$ in
standard form
\[
Q=\langle w_1\rangle \times \dots \times \langle w_k \rangle \times
G(Z).
\]
We call $(|Z|,k)$ the \emph{rank} of $H$. We use the left lexicographical
order on ranks of quasiparabolic subgroups: 
if $H$ and $K$ are 
quasiparabolic subgroups of ranks  $(|Z_H|,k_H)$ and 
$(|Z_K|, k_K)$, respectively,
then  $\rank(H) < \rank(K)$ if
 $(|Z_H|,k_H)$ precedes $(|Z_K|, k_K)$ in left
 lexicographical order.

The  centraliser $C_G(g)$ of an element $g\in G$ is a
typical example of a quasiparabolic subgroup
\cite{DK}. We shall see below
(Theorem \ref{thm:2bcentr}) that 
the centraliser of any set of elements of the group $G$ is
a quasiparabolic subgroup.

\begin{lem}\label{lem:qphom}
A quasiparabolic subgroup is a block-homogeneous subgroup and
consequently any intersection of quasiparabolic subgroups is again
block-homogeneous.
\end{lem}
\begin{proof}
Let $Q(w,Z)$ be a canonical quasiparabolic subgroup. Since
$w$ is a cyclically minimal root element it follows that $Q(w,Z)$ is 
 block-homogeneous. An application of Lemma \ref{lem:4.1} then implies
$Q(w,Z)^g$ is also block-homogeneous.
\end{proof}
We shall need the following lemma in Section \ref{sec:ht}.
\begin{lem}\label{lem:qpcont}
Let $Q_1=Q(u,Y)$ and $Q_2=Q(v,Z)$ be canonical quasiparabolic subgroups
in standard form
and let $g\in G$. If $Q_2^g\subseteq Q_1$, $g\in G(Z^\perp)$ and 
$\gd^r_Y(g)=1$ then $Q_2^g$ is a canonical quasiparabolic subgroup. 
\end{lem}
\begin{proof}
Let $u$ and $v$ have block decompositions $u=u_1\cdots u_k$ and 
$v=v_1\cdots v_l$, respectively. As $g\in G(Z^\perp)$ we have 
\[
Q_2^g= \la v_1^g\ra \times \cdots \times \la v_l^g\ra\times G(Z).
\]
Therefore, 
 for $j=1,\ldots ,l$, 
either $v_j^g=u_i$ for some $i=1,\dots, k$, 
or $v_j^g\in  G(Y)$. 
 If  $v_j^g=u_i$ then
$v_j^g$ is a cyclically minimal root element. 
If, on the other hand,  $v_j^g\in G(Y)$
then, from Lemma \ref{cor:conj}, there exist elements
$b,c,d$ and $e$ such that $g=b\circ c\circ d$, $v_j^g=d^{-1}\circ e\circ d$ and 
$e=v_j^b$ is a cyclically minimal root element. 
As $v_j^g\in G(Y)$  and $\gd_Y^r(g)=1$ we have 
$d=1$ and so $v_j^g=e$ and is a cyclically minimal root element. Therefore $Q_2^g$ is a
canonical quasiparabolic subgroup.
\end{proof}
\subsection{Intersections of Quasiparabolic Subgroups}

The main result of this section is the following
\begin{thm} \label{thm:intqpsub}
An intersection of quasiparabolic subgroups is a quasiparabolic
subgroup.
\end{thm}

We shall make use of the following results.
\begin{lem} \label{lem:4.2}
  Let $A=A_1\times \dots \times A_l$ and $B=B_1\times\dots \times
  B_k$, $A_i$, $B_j$, $i=1,\dots, l$, $j=1,\dots, k$ be block-homogeneous subgroups of $G$ and
  $C=A\cap B$. Then 
\[C=\prod\limits_{\begin{tiny}\begin{array}{c}
                                i=1,\dots, l; \\
                               j=1,\dots, k
                            \end{array}\end{tiny}
} (A_i\cap  B_j).\]
\end{lem}
\begin{proof}
  If $C=1$ then the result is straightforward. Assume then that $C\ne 1$,
  $w\in C$ and  $w\ne 1$ and let $w=w_1\dots w_t$ be the block decomposition
  of $w$. Since 
$C$ is a
  block-homogeneous subgroup, $w_i\in C$, $i=1,\dots, t$. 
As $w_i$ is a block element we have $w_i\in A_r$ and
  $w_i\in B_s$ and consequently $w_i$ lies in 
$\prod_{i,j}(A_i\cap B_j)$. As it is clear that 
$C\ge \prod_{i,j}(A_i\cap B_j)$
this proves the lemma.
\end{proof}

\begin{lem} \label{lem:conj}
  Let $Z\subseteq X$, $w\in G(Z)$, $g\in G$. Suppose that
  $u=g^{-1}wg$ is cyclically minimal and $\gcd_{\alpha(w)}^l(g)=1$,
  then $g$ and $w$ commute.
\end{lem}
\begin{proof}
  Let $g=d\circ g_1$, where $d=\gcd_{\alpha(w)^\perp}^l(g)$. If
  $g_1=1$ then $g\in C(w)$. Suppose $g_1\ne 1$. Then
  $\gcd_{\alpha(w)}^l(g_1)=1$ so we write $g_1=x \circ g_2$, where $x\in
  (X\cup X^{-1}) \bs (\alpha(w) \cup \alpha(w)^\perp)$
  and thus $u=g_2^{-1}x^{-1}w x g_2$ is written in geodesic form.
  This is a contradiction for $l(w)<l(u)$.
\end{proof}

\begin{lem}\label{lem:2quasi}
Let
\[ 
Q_1=\langle u_1 \rangle \times \dots \times \langle u_l \rangle
  \times G(Y) \textrm{ and } Q_2=\langle v_1 \rangle \times \dots \times \langle v_k
  \rangle
  \times G(Z)
\]
be canonical quasiparabolic subgroups in standard form and let $g\in G$
such that $\gd^l_Z(g)=1$.  
Write $g=d\circ h$, where $h=\gd^r_Y(g)$ and 
set $T=\a(d)^\perp$. 
Then, after reordering the $u_i$ and $v_j$ if necessary, there exist
$m, s, t$ such that  
\begin{equation} \label{eq:intqp}
Q_1\cap Q_2^g=\left(\prod_{i=1}^s \la v_i\ra\times 
\prod_{i=s+1}^t\la v_i\ra \times \prod_{j=s+1}^m\la u_i\ra
\times G(Y\cap Z\cap T)\right)^g
\end{equation}
and
\be[(i)]
\item $\la u_i\ra =\la v_i\ra^g$, for $i=1,\ldots ,s$;
\item\label{item:2quasi2} $\la v_i\ra^g\subseteq G(Y)$, for $i=s+1,\ldots ,t$; and 
\item $\la u_i\ra \subseteq G(Z)$, for $i=s+1,\ldots ,m$.
\ee
\end{lem}
\begin{proof}
As $Q_1\cap Q_2^g=(Q_1\cap Q_2^d)^h$ we may assume that $h=1$ and $d=g$,
so $\gd^r_Y(g)=1$. As $Q_1$ and $Q_2^g$ are block-homogeneous 
we may apply Lemma \ref{lem:4.2} to compute their intersection.
Therefore we consider the various possible intersections of factors of $Q_1$ and
$Q_2^g$. 
\be[(i)]
\item If $\langle u_i\rangle \cap {\langle v_j\rangle}^g\ne 1$ then, as 
$u_i$ and $v_j$ are root elements, 
$\langle u_i \rangle =\langle v_j\rangle ^g$. 
Suppose that this is the case for $u_1,\dots u_s$ and
$v_1,\dots, v_s$ and that $\la u_i\ra \cap \la v_j\ra ^g=1$, if $i>s$ or $j>s$.
\item If $\langle v_j\rangle ^g\cap G(Y) \ne 1$ then, since $v_j$ is cyclically minimal, $\langle
    v_j\rangle^g \subset G(Y)$. This cannot happen if $j\le s$ so suppose it 
is the case for
    $v_{s+1}, \dots, v_t$, and that $\la v_j \ra^g\cap G(Y)=1$, for $j>t$.
\item If $\langle u_i\rangle \cap {G(Z)}^g \ne 1$ then
    $u_i=w^g$, $w\in G(Z)$ and by Lemma \ref{lem:conj}, $w$ and $g$
    commute so does $u_i=w=u_i^g$. This cannot happen if $i\le s$ so 
suppose that it's the case for
    $u_{s+1},\dots,u_{m}$, and not for $i>m$.
\item Finally, using Lemma \ref{lem:formpar} and the assumption that
$\gd^r_Y(g)=\gd^l_Z(g)=1$,  we have
$G(Y)\cap {G(Z)}^g={G(Y\cap Z\cap T)}={G(Y\cap Z\cap T)}^g$,
where $T=\alpha(g)^\perp$.
\ee
Combining these intersections \eqref{eq:intqp} follows from Lemma \ref{lem:4.2}.
\end{proof}

\begin{cor}\label{cor:qprank}
Let $H_1$ and $H_2$ be quasiparabolic subgroups of $G$  then $H_1\cap H_2$ is
quasiparabolic and  {\rm
$\rank(H_1\cap H_2) \le \min\{\rank(H_1), \rank(H_2)\}$}.
\end{cor}
\begin{proof}
Let $H_1=Q_1^f$ and $H_2=Q_2^g$, where $Q_1=Q(u,Y)$ and $Q_2=Q(v,Z)$ are quasiparabolic
subgroups in standard form, as in Lemma \ref{lem:2quasi}.
As in the proof of Proposition \ref{prop:intpar} we may assume that 
$f=1$ and $\gd^l_Z(g)=1$ and so Lemma \ref{lem:2quasi} implies 
$H_1\cap H_2$ is quasiparabolic. If $\rank(H_1\cap H_2)\ge \rank(H_1)$
then $|Y|\le |Y\cap Z\cap T|$ so $Y\subseteq Z\cap T$. In this case 
\eqref{item:2quasi2} of Lemma \ref{lem:2quasi} cannot occur.
Therefore, in the notation of Lemma \ref{lem:2quasi}, 
$\rank(H_1\cap H_2)=s+m$. If $\rank(H_1\cap H_2)\ge \rank(H_1)$ then
$s+m\ge l$ which implies $m=l-s$ and so $u_i\in G(Z)^g$, for $i=s+1,\ldots ,l$. 
As $u_i=v_i^g$, for $i=1,\ldots, s$ it follows that $H_1\subseteq H_2$.
\end{proof}
\begin{proof}[Proof of Theorem \ref{thm:intqpsub}.]
  Given Corollary \ref{cor:qprank} the intersection of 
an infinite collection of quasiparabolic
subgroups is equal to  the intersection of a finite sub-collection.
From Corollary \ref{cor:qprank} again such an intersection is quasiparabolic
and the result follows.
\end{proof}

\subsection{A Criterion for a Subgroup to be a Centraliser} \label{sec:2bcentr}

\begin{thm} \label{thm:2bcentr}
A subgroup $H$ of $G$ is a centraliser if and only if the two
following conditions hold.
\begin{enumerate}
    \item \label{cond1} $H$ is conjugate to some canonical
    quasiparabolic subgroup $Q$.
    \item \label{cond2} If $Q$ is written in standard form
    $$
    Q=\langle w_1\rangle \times \dots \times \langle w_k \rangle \times
    G(Y),
    $$
where
    $w=w_1\dots w_k$ is the block decomposition of a cyclically
    minimal element $w$, $w_i$ is a root element and  
$|\alpha(w_i)|\ge 2$, $i=1,\dots, k$, then 
    \[
Y\in \CS(\G)\textrm{ and }Y\in \CS(\G_w)\textrm{ 
where }\G_w=\G(\alpha(w)^\perp\bs\a(w)).\]
\end{enumerate}
\end{thm}
\begin{proof}
  Let $H=C(u_1,\dots, u_l)$. Then $H=\bigcap\limits_{i=1}^k C(u_i)$ and we
  may assume that each $u_i$ is a block root element. Since $C(u_i)$ is a
  quasiparabolic subgroup, then by Theorem \ref{thm:intqpsub}, $H$
  is also a quasiparabolic subgroup and is conjugate to a canonical
  quasiparabolic subgroup $Q=\langle w_1\rangle \times \dots \times \langle w_k \rangle \times
    G(Y)$ written in standard form. Thus condition \ref{cond1} is
    satisfied. 

Then $H=Q^g$ and, after  conjugating all the $u_i$'s by
    $g^{-1}$ we have  a centraliser $H^{g^{-1}}=Q$. Thus we may assume that 
$H=Q$.
Let $w=w_1\cdots w_k$, set $Z=\alpha(w)^\perp\bs\a(w)$ and 
$T=\bigcup\limits_{i=1}^l \alpha (u_i)$.
As $w$ has block decomposition $w=w_1\cdots w_k$ we have 
$C(w)=\langle w_1\rangle \times \dots \times
\langle w_k \rangle \times G(Z)$.
For all $y\in Y$ we have $y\in C(u_i)$ so 
and thus $y\in C(\alpha(u_i))$ and  
  $Y\subseteq T^\perp$.
Conversely if  $y\in T^\perp$ then $y\in C(u_i)$ 
so $y\in Q$ and, by definition of standard form, $y\in Y$. Therefore $Y=T^\perp$. 
It follows that  $Y\in \CS(\G)$
and since by Lemma \ref{lem:stdform} we have $Y\cap \a(w)=\nul$ we also
have $Y\subseteq Z$.

It remains to prove that $Y\in \CS(\G_w)=\CS(\G(Z))$. 
Set $W=\alpha(w)$. We show that $T\cup Z\subseteq W\cup Z$. Take
$t\in T=\bigcup\limits_{i=1}^l \alpha (u_i)$, $t\notin W$ and
suppose that $t\in \alpha(u_m)$. Since $w \in C(u_i)$, we have $u_m
\in C(w)=\langle w_1\rangle \times \dots \times \langle w_k \rangle
\times G(Z)$. Now $u_m$ is a root block element and 
$C(w)$ is a block-homogeneous subgroup so if 
$ u_m = w_j^{\pm 1}$ for some $j$ then $t \in W=\alpha(w)$, contrary to 
the choice of $t$. Therefore $u_m \in G(Z)$, so $t\in Z$  and $T \cup
Z\subseteq W\cup Z$, as claimed.

Assume now that $Y\notin \CS(\G(Z))$. In this case there exists an element
$z\in Z\bs Y$ such that $z\in \cl_Z(Y)$. Since $z\notin Y=T^\perp$,
there exists $u_m$ such that $[u_m,z]\ne 1$ and so there exists $t\in \alpha(u_m)$
such
that $[t, z]\ne 1$. As $[z,W]=1$, we have $t\notin W$ and since
$W\cup Z\supseteq T\cup Z$, we get $t\in Z$. This together with
$t\in \alpha(u_m)\subseteq Y^\perp$ implies that $t\in \cO^Z(Y)$. Since
$[z,t]\ne 1$, we obtain $z\notin \cl_Z(Y)$, in contradiction to 
the choice of $z$. Hence $\cl_Z(Y)=Y$ and $Y\in \CS(\G(Z))$.

Conversely, let $Q=\langle w_1\rangle \times \dots \times \langle
w_k \rangle \times G(Y)$ be a canonical quasiparabolic subgroup
written in the standard form, $Y\in \CS(X)$ and $Y \in \CS(\G(Z))$,
where $Z=\a(w)^\perp\bs\a(w)$.
We shall prove that $Q=C(w, z_1,\dots, z_l)$, where $z_1,\dots, z_l$ are
some elements of $Z$.  If $Y=Z$ then $Q=C(w)$. If $Y\subsetneq Z$
then, since $Y=\cl_Z(Y)$, there exist $z_1,u\in Z$ so that $z_1\in
\cO^Z(Y)$ and $[z_1,u]\ne 1$. In which case $C(w,z_1)=\langle
w_1\rangle \times \dots \times \langle w_k \rangle \times G(Y_1)$,
$Y\subseteq Y_1 \subsetneq Z$ (the latter inclusion is strict for
$u\notin Y_1$). If $Y_1=Y$ then $Q=C(w, z_1)$, otherwise iterating
the procedure above, the statement follows.
\end{proof}

A centraliser which is equal to a canonical quasiparabolic subgroup
is called a {\em canonical quasiparabolic centraliser}.
\section{Height of the Centraliser Lattice}\label{sec:ht}
In this section we will give a new shorter proof of the main theorem
of \cite{DKR2}.

\begin{thm}\label{thm:height}
Let $G=G(\G)$ be a free partially commutative group, let $\fC(G)$ be
its centraliser lattice and let $L=\CS(\G)$ be the lattice of closed
sets of $\G$. Then the height $h(\fC(G))=m$ equals the height $h(L)$
of the lattice of closed sets $L$.
\end{thm}

In order to prove this theorem we introduce some notation for the various parts of canonical
quasiparabolic subgroups.
\begin{defn}
Let $Q=\la w_1\ra\times \cdots \times \la w_k\ra\times G(Z) $ be 
a quasiparabolic
subgroup in standard form. Define the {\em block set} of $Q$ to be  $\cB(Q)=
\{w_1 \ldots, w_k\}$ and the 
{\em parabolic part} of $Q$ to be  $\cP(Q)=G(Z)$. 
Let $Q^\prime $ be a quasiparabolic subgroup with block set
$\la v_1\ra \times \cdots \times \la v_l\ra$ and parabolic part $G(Y)$.
Define the {\em block difference} of $Q$ and $Q^\prime$ to be
$b(Q,Q^\prime)=|\cB(Q)\bs \cB(Q^\prime)|$, that is the number of blocks occurring in
the block set of $Q$ but not $Q^\prime$. Define the {\em parabolic difference}
of $Q$ and $Q^\prime$ to be $p(Q,Q^\prime)=|Z\bs Y|$.
\end{defn}%

The following lemma is the key to the proof of the theorem above.
\begin{lem}\label{lem:cpad}
Let $C$ and $D$ be canonical quasiparabolic centralisers such that
$C>D$ and $b=b(D,C)>0$.
Then $p(C,D)>0$ and there exists a strictly descending chain 
of canonical parabolic centralisers 
\begin{equation}\label{eq:keychain}
\cP(C)>C_b>\cdots >C_1>\cP(D)
\end{equation}
of length $b+1$.
\end{lem}
\begin{proof}
Let
 $C$ and $D$ have parabolic parts
$\cP(C)=G(Y)$ and $\cP(D)=G(Z)$, for closed subsets $Y$ and $Z$ in $\CS(\G)$.
Let the block sets of $C$ and $D$ be
$\cB(C)=\{u_1,\ldots ,u_k\}$ and $\cB(D)=\{v_1,\ldots , v_l\}$. Fix
$i$ with $1\le i\le l$. As 
$D<C$, either $\la v_i\ra=\la u_j\ra$, for some $j$, or $\la v_i\ra\subseteq
G(Y)$. As $b(D,C)>0$ there exists $i$ such that $\la v_i\ra\subseteq G(Y)$.
Moreover, for such $i$,  we have $\a(v_i)\subseteq Y\bs Z$, so 
$p(C,D)>0$. 

Assume that, after relabelling if necessary, 
$\la v_i\ra =\la u_i\ra$, for $i=1,\ldots ,s$, and that
$\la v_{s+1}\ra,\ldots ,\la v_{l}\ra\subseteq G(Y)$, so $b=l-s$.  Choose 
$t_i\in \a(v_{s+i})$ and let $Y_i=\cl(Z\cup\{t_1,\ldots ,t_i\})$, 
for $i=1,\ldots ,l-s=b$.
  Let $C_i=G(Y_i)=C_G(Y_i^\perp)$, so $C_i$ is a canonical
parabolic centraliser. We claim that the chain \eqref{eq:keychain} is 
strictly descending. To begin with, as $t_1\in Y_1\bs Z$ we have $G(Z)<C_1$.
Now fix $i$ and $n$ such that $1\le i < n \le b$. If $a\in \a(v_{s+n})$ then
$a\in Z^\perp$ and $a\in \a(v_{s+j})^\perp\subseteq t_j^\perp$, for $1\le j<n$,
by definition of the standard form of quasiparabolic subgroups. Hence
$a\in (Z\cup\{t_1,\ldots ,t_i\})^\perp$. Thus $[a,b]=1$, for all
$b\in Y_i$. This holds for all $a\in \a(v_{s+n})$ so 
$Y_i\subseteq \a(v_{s+n})^\perp$. As $v_{s+n}$ is a block of length
at least $2$ we have $\a(v_{s+n})\cap \a(v_{s+n})^\perp=\emptyset$, so
$Y_i\cap \a(v_{s+n})=\emptyset$. Hence $t_n\notin Y_i$ and it follows
that $C_i<C_{i+1}$, $i=1,\ldots ,b-1$. Now choose $c\in \a(v_{s+1})$ 
such that $[c,t_1]\neq 1$. Then $c\in Y$, as $\a(v_{s+1})\subseteq Y$,
however $c\notin Y_b$, since $t_1\in Z^\perp \cap t_1^\perp\cap \cdots
\cap t_b^\perp=Y_b^\perp$ and $Y_b=Y_b^{\perp\perp}$. As $D<C$ we have 
$Z\subseteq Y$ so 
$C_b=G(Y_b)<G(Y)$.
\end{proof}
We can use this lemma to prove the following about chains of 
canonical quasiparabolic subgroups.
\begin{lem}\label{lem:qparac}
Let $C_0>\cdots >C_d$ be a strictly descending chain of 
canonical quasiparabolic centralisers such that $C_0$ and
$C_d$ are canonical parabolic centralisers. Then there exists
a strictly descending chain $C_0>P_1>\cdots >P_{d-1}>C_d$, of
canonical parabolic centralisers.
\end{lem}
\begin{proof}
First we divide the given centraliser chain 
into types depending on
block differences. Then we replace the chain with a chain of 
canonical parabolic centralisers, using Lemma \ref{lem:cpad}. 
A simple counting argument shows that the new chain has 
length at least as great as the old one. 
In detail let $I=\{0,\ldots ,d-1\}$ and 
\begin{align*}
I_+ &=\{i\in I: b(C_{i+1},C_i)>0\},\\
I_0 &=\{i\in I: b(C_{i+1},C_i)=0\textrm{ and } p(C_i,C_{i+1})>0\} 
\textrm{ and }\\
I_- &= \{i\in I: b(C_{i+1},C_i)= p(C_i,C_{i+1})=0\}.
\end{align*}
Then $I=I_+\sqcup I_0\sqcup I_-$. For $i\in I_+$ let $\D_i$ be the   strictly
descending chain of canonical parabolic centralisers of length 
$b(C_{i+1},C_i)+1$ from $\cP(C_i)$ to $\cP(C_{i+1})$, constructed in 
Lemma \ref{lem:cpad}. For $i\in I_0$ let $\D_i$ be the length one chain
$\cP(C_i)>\cP(C_{i+1})$ and for $i\in I_-$ let $\D_i$ be the length
zero chain $\cP(C_i)=\cP(C_{i+1})$. This associates 
a chain $\D_i$ of canonical parabolic centralisers to each $i\in I$ 
and we write $l_i$ for the length of $\D_i$.
If $\D_i=P_0>\cdots >P_{l_i}$ and 
$\D_{i+1}=P^\prime_0>\cdots >P^\prime_{l_{i+1}}$ then 
by definition $P_{l_i}=P^\prime_0$, for $i=1,\ldots ,d-1$. We may
therefore concatenate $\D_i$ and $\D_{i+1}$ to give a chain
of canonical parabolic centralisers
\[P_0> \cdots >P_{l_i}=P^\prime_0>\cdots >P^\prime_{l_{i+1}}\]
of length $l_i+l_{i+1}$. Concatenating $\D_1, \ldots ,\D_{d-1}$ in
this way we obtain a strictly descending chain of canonical parabolic centralisers
of length $l=\sum_{i=0}^{d-1} l_i$.
Moreover 
\[l=\sum_{i\in I_+}b(C_{i+1},C_i)+|I_+|+|I_0|,\]
since $l_i=b(C_{i+1},C_i)+1$, for all $i\in I_+$, $l_i=1$, for all $i\in I_0$ and $l_i=0$, for all $l_i\in I_-$. As $|I|=d$ we have now
\[
l-d =\sum_{i\in I_+}b(C_{i+1},C_i)-|I_-|.\]
To complete the argument we shall show that
\[
\sum_{i\in I_+}b(C_{i+1},C_i)=|\cup_{i=0}^d\cB(C_i)|\ge |I_-|.
\]
As $\cB(C_0)=\emptyset$ we have $b(C_1,C_0)=|\cB(C_0)\cup\cB(C_1)|$.
Assume inductively that 
\[
\sum_{i=0}^k b(C_{i+1},C_i)=|\cup_{i=0}^{k+1}\cB(C_i)|, 
\]
for
some $k\ge 0$. Then 
\[
\sum_{i=0}^{k+1} b(C_{i+1},C_i)=
|\cup_{i=0}^{k+1}\cB(C_i)|+|\cB(C_{k+2})\bs \cB(C_{k+1})|.
\]
Moreover, if $w\in \cB(C_{k+2})\bs \cB(C_{k+1})$ then 
$w\in \cP(C_j)$, for all $j\le k+1$, so $w\notin \cB(C_j)$, for
$j=0,\ldots ,k+1$. Hence 
\[
\cB(C_{k+2})\bs \cB(C_{k+1})=\cB(C_{k+2})\bs \cup_{i=0}^{k+1}\cB(C_i)\]
and it follows that 
\[
\sum_{i=0}^{k+1} b(C_{i+1},C_i)=|\cup_{i=0}^{k+2}\cB(C_i)|. 
\]
As $b(C_{i+1},C_i)=0$ if $i\notin I_+$ it follows that  
\[
\sum_{i\in I_+}b(C_{i+1},C_i)=|\cup_{i=0}^d\cB(C_i)|,
\]
 as required. If $i\in I_-$ then $b(C_{i+1},C_i)=p(C_i,C_{i+1})=0$, so
$b(C_i,C_{i+1})>0$. Therefore there is at least one element $w\in 
\cB(C_i)\bs \cB(C_{i+1})$. It follows that $w\notin \cB(C_j)$, for
all $j\ge i+1$ and so $I_-\le |\cup_{i=0}^d\cB(C_i)|$. Therefore
$l-d\ge 0$ and the proof is complete.
\end{proof}
\begin{proof}[Proof of Theorem \ref{thm:height}]
Let 
\begin{equation*}
G=C_0>\cdots > C_d=Z(G) 
\end{equation*}
be a maximal descending chain of
centralisers of $G$. By Theorem \ref{thm:2bcentr}, each of the
$C_i$'s is a quasiparabolic subgroup. 
If each 
$C_i$ is canonical 
then, since $G$ and $Z(G)$ are both canonical parabolic centralisers
the result follows from Lemma \ref{lem:qparac}.

Suppose  that now $C_1,\dots, C_s$ are canonical quasiparabolic and
$C_{s+1}$ is not: say 
$C_{s+1}=Q^{g}$, where $Q$ is a canonical
quasiparabolic subgroup.
Let $C_s=Q(u,Y)$ and $Q=Q(v,Z)$ both in standard form. 
Write  $g=f\circ h$, where $f=\gd_{Z^\perp}^l(g)$ and let
$f=e\circ d$, where $d=\gd^r_Y(f)$, so $d \in G(Y\cap Z^\perp)$. Then
$G(Z)^h=G(Z)^{dh}= G(Z)^g\subseteq G(Y)$ and $\a(h)\subseteq Y$, from Corollary 
\ref{cor:parin}. Hence $\a(d\circ h)\subseteq Y$ which implies that  
$\a(d\circ h)\subseteq \cP(C_s)\subseteq \cdots
\subseteq \cP(C_0)$. It follows that $C_r^{dh}=C_r$, for $r=0,\ldots ,s$.
Therefore conjugating $C_0>C_1>\cdots >C_d$ by
$(dh)^{-1}$ we obtain a chain in which $C_0,\ldots, C_{s}$ are unchanged
and $C_{s+1}=Q^e=\la v_1\ra^e\times \la v_l \ra^e\times G(Z)$,
with $\gd^l_Z(e)=\gd^r_Y(e)=1$, $e\in G(Z^\perp)$.
As Lemma \ref{lem:qpcont} implies that  $Q^e$ is a
canonical quasiparabolic subgroup we now have a chain in which
$C_0,\ldots, C_{s+1}$ are canonical quasiparabolic. Continuing this
way we eventually obtain a chain, of length $d$, 
of canonical quasiparabolic centralisers 
to which the first part of the proof may be applied. 
\end{proof}

\end{document}